# ON A CONJECTURED INEQUALITY IN CONVEX ANALYSIS FOR THE CASE OF THE UNIT BALL OF $l_p^n$, $1 \leq p \leq \infty$

## D. KARAYANNAKIS

**Abstract**

If $B_p^n$ is the unit $p$-ball of the Euclidean $(R^n, \langle \ \rangle)$, $|B_p^n|$ its volume and $q$ the Hölder conjugate exponent of $p$, we prove that $\dfrac{\int_{B_p^n} \int_{B_q^n} \langle x, y \rangle^2 \, dxdy}{|B_p^n||B_q^n|} \leq \dfrac{n}{(n+2)^2}$. This confirms one of Kuperberg's recent conjectures for centrally symmetric convex bodies $K$ of $R^n$ in the special case $K = B_p^n$, by the exclusive use of basic classical analysis.

## Part I: Notation, goals and main strategy

One of the key quantities in convex analysis is the ratio $f(K) = \dfrac{1}{|K||K^o|} \int_{K \times K^o} \langle x, y \rangle^2 \, dxdy$,

where $K$ is a (centrally symmetric) convex compact set of the Euclidean inner-product space $(R^n, \langle \ \rangle)$ with non-empty interior (or simply "body") and

$K^o = \{y \in R^n : |\langle x, y \rangle| \leq 1\}, \forall x \in K\}$ is the so-called polar set of $K$. This quantity is tightly connected with various key notions of convex geometric analysis like the isotropic constant, the hyperplane conjecture and the Blashke-Santaló inequality (for a survey see e.g [7] or [2]). Recently G. Kuperberg in [6], among others, conjectured in a probabilistic formulation that for any symmetric body $K$ we have that $f(K) \leq f(B_2^n) = \dfrac{n}{(n+2)^2}$.

A proof of this conjecture for the special case $K = B_p^n$ was given very recently by D.A. Gutiérrez in [3] by use of polygamma functions and pertinent convexity theorems. In the present work we also confirm the above inequality by appealing only to the basic features and properties of the gamma function (as tabulated e.g. in [2]), along with aspects of fundamental classical analysis.

To achieve this goal we start in Lem.II.1 with a representation of $f(B_p^n) = f(n, p)$ by use of the gamma function and then, based on a suitable gamma function ratio formula (Lem.II.2) established in [5] (and already used much earlier for other purposes and only through a formal reasoning as e.g. in [4]), we obtain in Prop.II.1 a representation of $f(n, p)$ free of its gamma symbolism. Then proving in Prop.II.2 that $\dfrac{df(n, p)}{dp} \geq 0$, for all $n$ and $1 \leq p \leq 2$ we arrive, due to an obvious algebraic symmetry, at the announced result for $1 \leq p \leq \infty$. The method also furnishes manageable lower and upper bounds of $f(B_p^n)$ for $p_1 \leq p \leq p_2$ when $(p_1, p_2) \subset [1, 2]$ or $(p_1, p_2) \subset [2, \infty]$ (Corol.II.1 & Corol.II.2 respectively).

**Part II: A representation of $f(n, p)$ and its monotonicity**

**Lemma II.1**

(a) For $1 < p < \infty$, $f(n, p) = f(n, q) = \dfrac{\int_{B_p^n}\int_{B_q^n} \langle x, y \rangle^2 \, dxdy}{|B_p^n||B_q^n|} = n \dfrac{\Gamma(\frac{3}{p})\Gamma(\frac{3}{q})\Gamma(1+\frac{n}{p})\Gamma(1+\frac{n}{q})}{\Gamma(\frac{1}{p})\Gamma(\frac{1}{q})\Gamma(1+\frac{n+2}{p})\Gamma(1+\frac{n+2}{q})}$.

(b) $f(n,1) = f(1,\infty) = \dfrac{2n}{3(n+1)(n+2)}$.

Proof:

It is trivial that $f(B_p^1) = \dfrac{1}{9}$ and so both cases are directly verifiable Let now $n \geq 2$.

The formula for the volume $|B_p^n|$ is well-known, dating back to Dirichlet by use of multiple integration (see e.g. [1]), namely $|B_p^n| = \dfrac{[2\Gamma(1+\frac{1}{p})]^n}{\Gamma(1+\frac{n}{p})}$, (1)

(and note also that $|B_1^n| = \dfrac{2^n}{n!}$ and $|B_\infty^n| = 2^n$), so we concentrate on the double integral of $f(n, p)$ separately in each case.

For (a), under the evident abuse of notation, setting $\phi(p) = \int_{B_p^n} x_1^2 dx$ and due to the symmetry with respect to the basic planes, we have that $\iint_{B_p^n \times B_q^n} \langle x, y \rangle^2 dxdy = n\phi(p)\phi(q)$. On the other hand

$$\phi(p) = 2\int_0^1 t^2 |B_p^{n-1}| (1-t^p)^{\frac{n-1}{p}} dt = \dfrac{2}{p}|B_p^{n-1}| \int_0^1 \dfrac{s^{\frac{2}{p}}(1-s)^{\frac{n-1}{p}} s^{\frac{1}{p}-1}}{s^{\frac{3}{p}-1}} ds,$$

by the change of variable $s = t^p$.

Using Euler's first integral (or beta function representation-see e.g. [1]) we have that

$$\phi(p) = \dfrac{2}{p}|B_p^{n-1}| \dfrac{\Gamma(\frac{3}{p})\Gamma(1+\frac{n-1}{p})}{\Gamma(1+\frac{n+2}{p})}. \quad (2)$$

Let (1)', (2)' the corresponding formulas for $|B_q^n|$ and $\phi(q)$ respectively. Combining (1), (2), (1)' and (2)' and using the identity $\Gamma(x+1) = x\Gamma(x)$ we obtain the announced result.

For (b) we have $\phi(\infty) = 2\int_0^1 t^2 |B_\infty^{n-1}| dt = \dfrac{2^n}{3}$ and similarly, along the lines of part (a),

$\phi(1) = \dfrac{2^{n+1}}{(n+2)!}$. Thus by substitution in $f(n,1)$ we obtain the announced result. $\square$

In [5] we have proved the following:

**Lemma II.2**

For $x > 0$, $a < 1$, $x + a \notin \mathbb{Z}_0^-$, $\quad \Gamma(1-a)\dfrac{\Gamma(x+a)}{\Gamma(x)} = \displaystyle\prod_{k=1}^{\infty}\dfrac{k(k+x-1)}{(k-a)(k+x+a-1)}.\square$

**Remark II.1**

Naturally when in (a) $p \to 1^+$ we must obtain (b). For this it will suffice to show that

$$\lim_{q\to\infty}\left[\dfrac{\Gamma(\frac{3}{q})\,\Gamma(\frac{n}{q})}{\Gamma(\frac{1}{q})\,\Gamma(\frac{n+2}{q})}\right] = \dfrac{n+2}{3} \quad \text{which can be verified directly using induction on } n \text{ together}$$

with Lem.II.2 for $x = \dfrac{n}{q}$, $x = \dfrac{n+2}{q}$ and $a = \dfrac{1}{q}$. $\square$

Combining now the two above lemmas we obtain :

**Proposition II.1**

For $1 \leq p \leq \infty$, $f(n, p) = \dfrac{n}{9}\displaystyle\prod_{k=1}^{\infty}\dfrac{(k^2+k+\frac{1}{pq})\,[k^2+(n+2)k+\frac{(n+2)^2}{pq}]}{(k^2+3k+\frac{9}{pq})\,(k^2+nk+\frac{n^2}{pq})}$

Proof:

For $p = 1$ ($q = \infty$) or $p = \infty$ ($q = 1$) the above infinite product, after telescoping, reduces to $\dfrac{6}{(n+1)(n+2)}$ and thus we obtain the formula of Lem.II.1(b).

For $1 < p < \infty$ we only need to apply Lem.II.2 in the numerator (resp. denominator) of the infinite product of Lem.II.1(a) for $x = \dfrac{2}{p}, a = \dfrac{1}{p}$ and $x = \dfrac{2}{q}, a = \dfrac{1}{q}$ (resp. $x = \dfrac{n+1}{p}, a = \dfrac{1}{p}$ and $x = \dfrac{n+1}{q}, a = \dfrac{1}{q}$) and then simplify inside the infinite product. $\square$

Exploiting the above representation we will show that, for any fixed $n$, $f(n,p)$ is (strictly for $n \neq 1$) increasing over $[1, 2]$.

**Proposition II.2**

(a) For $1 < p \leq 2$, $\dfrac{df(n,p)}{dp} \geq 0$, $\forall n \in N$, with equality iff $p = 2$ when $n \neq 1$.

(b) For all $n$ and $1 \leq p \leq 2$, $\max f(n,p) = \dfrac{n}{(n+2)^2}$.

Proof:

(a) Since $f(1,p) = \dfrac{1}{9}$ we can consider $n \geq 2$. Let $t = \dfrac{1}{pq}$. Then $t > 0$ (and $t \leq \dfrac{1}{4}$). Since

$\dfrac{dt}{dp} = \dfrac{2-p}{p^3} \geq 0$ (with equality iff $p = 2$) it will suffice to show $f' = \dfrac{df}{dt} \geq 0$ over $(0, \dfrac{1}{4}]$, for all $n \geq 2$ (under the evident abuse of notation).

Let $g_k(n,t) = k^2 + nk + n^2 t$.

Then, since again under the evident abuse of notation $g_k'(n) = n^2$, it is straightforward through logarithmic differentiation to see that

$$\text{sgn } f' = \text{sgn} \sum_{k=1}^{\infty} \left[ \frac{n^2 g_k(3) - 9}{g_k(1) g_k(3)} + \frac{(n+2)^2 g_k(n) - n^2 g_k(n+2)}{g_k(n) g_k(n+2)} \right].$$

Now we observe that for any fixed $k$ and $0 < t \leq \frac{1}{4}$ (and in fact for any $t > 0$)

$$n^2 g_k(3) g_k(n+2)[g_k(n) - g_k(1)] + g_k(n)[(n+2)^2 g_k(1) g_k(3) - 9 g_k(n+2)] > 0, \quad (3)$$

based on the facts that the difference in the first bracket of (3) is evidently positive and the one in the second bracket is also positive by simple induction on $n$. □

(b) In face of (a), Lem.II.1(b), Rem.II.1 and Prop.II.1 we can check by telescoping that

$$f(n,1) = \frac{2n}{3(n+1)(n+2)} \leq f(n, \tfrac{1}{4}) = \frac{n}{9} \left[ \prod_{k=1}^{\infty} \left( \frac{2k+1}{2k+3} \right) \left( \frac{2k+n+2}{2k+n} \right) \right]^2 = \frac{n}{(n+2)^2}. \square$$

### Remark II.2

By Prop.II.1, $f(n,p) = f(n,q)$, and since $1 \leq p \leq 2$ iff $2 \leq q \leq \infty$ we see that Prop.II.2 holds for $1 \leq p \leq \infty$. In addition the proof of Prop.II.2(a) shows that $f(n,p)$ is (strictly for $n \neq 1$) decreasing over $(2, \infty)$ and thus we are led to the following inequalities:

### Corollary II.1

Let $n \geq 2$ and $1 \leq r < s \leq 2$. Then $\prod_{k=1}^{\infty} \frac{g_k(1,R) g_k(n+2,R)}{g_k(3,R) g_k(n,R)} < \prod_{k=1}^{\infty} \frac{g_k(1,S) g_k(n+2,S)}{g_k(3,S) g_k(n,S)}$ where $R = \frac{r-1}{r^2}$ and $S = \frac{s-1}{s^2}$.

### Corollary II.2

Let $n \geq 2$ and $2 \leq r < s \leq \infty$. Then the inequality in Corol.II.1 is reversed.

DEPARTMENT OF SCIENCES , TEI OF CRETE, 71004 HERAKLION, GREECE
*E-mail address*: dkar@stef.teiher.gr